\documentclass{amsart} 
\usepackage{enumerate} 
\usepackage{upref}
\usepackage{verbatim} 
\usepackage{color}
\usepackage{graphicx}
\usepackage[textsize=tiny]{todonotes}

\newcommand{\sign}{\mathop{\rm sign}}
\newcommand*{\mailto}[1]{\href{mailto:#1}{\nolinkurl{#1}}}

\usepackage[latin1]{inputenc}
\usepackage{amssymb, amsmath, amsfonts, amsthm}
\usepackage[all]{xy}

\usepackage{marginnote}

\DeclareMathOperator{\id}{Id}

\DeclareMathOperator{\supp}{supp}
\newcommand{\dott}{\, \cdot\,}

\newcommand{\Gr}{G}

\newcommand{\D}{\ensuremath{\mathcal{D}}}

\newcommand{\F}{\ensuremath{\mathcal{F}}}

\newcommand{\R}{\ensuremath{\mathcal{R}}}

\newcommand{\Real}{\mathbb R}

\newcommand{\muac}{\mu_{\text{\rm ac}}}

\DeclareMathOperator{\sgn}{sgn}

\newtheorem{theorem}{Theorem}[section]

\newtheorem{definition}[theorem]{Definition}

\numberwithin{equation}{section}

\allowdisplaybreaks

\begin{document}

\title[Peakon-antipeakon example]{The  general peakon-antipeakon solution for the Camassa--Holm equation}

\author[K. Grunert]{Katrin Grunert}
\address{Department of Mathematical Sciences\\ Norwegian University of Science and Technology\\ NO-7491 Trondheim\\ Norway}
\email{\mailto{katring@math.ntnu.no}}
\urladdr{\url{http://www.math.ntnu.no/~katring/}}

\author[H. Holden]{Helge Holden}
\address{Department of Mathematical Sciences\\
  Norwegian University of Science and Technology\\
  NO-7491 Trondheim\\ Norway}
\email{\mailto{holden@math.ntnu.no}}
\urladdr{\url{http://www.math.ntnu.no/~holden/}}

\thanks{Research supported in part by the
  Research Council of Norway  and by the Austrian Science Fund (FWF) under Grant No.~J3147.}  

\begin{abstract}
We compute explicitly the peakon-antipeakon solution of the Camassa--Holm equation $u_t-u_{txx}+3uu_x-2u_xu_{xx}-uu_{xxx}=0$ in the non-symmetric and $\alpha$-dissipative case. The solution 
experiences wave breaking in finite time, and the explicit solution illuminates the interplay between the various variables.
 \end{abstract}
\maketitle

\section{Introduction}
The Camassa--Holm (CH) equation
\begin{equation}\label{eq:CH}
 u_t-u_{txx}+3uu_x-2u_xu_{xx}-uu_{xxx}=0
 \end{equation}
 was first studied in the context of water waves in the seminal papers \cite{CH:93,CHH:94}. It possesses many interesting properties, including complete integrability and soliton-like solutions that interact in a manner similar to the solitons of the KdV equation.  In the context of the CH equation, solitons go by the name of peakons and antipeakons, and they are the topic of interest in the current paper. The peakons, that are stable solutions \cite{consstrauss},  are considerably more challenging than the KdV solitons as the peakons experience wave breaking in finite time and become singular.
 
The explicit example of the peakon-antipeakon solution for the Camassa--Holm equation has been a constant source of 
inspiration and intuition for the analysis of the solution of the general Cauchy problem.  As in a laboratory, one can test one's intuition on this particular solution that encodes most of the intricacies of the Cauchy problem. The key question is to analyze the behavior of the solution $u$ near wave breaking where 
$u_x(x_0,t)\to-\infty$ as $t\to t_0$, yet the $H^1$ norm remains finite \cite{cons_esc1}. Multipeakons can even be used for numerical computations for the general Cauchy problem, see \cite{HolRey:2006}.

Multipeakons appear as linear combinations of single peakons of the form
\begin{equation*}
 u(x,t)=\sum_{i=1}^np_i(t) e^{-\vert x-q_i(t)\vert}.
 \end{equation*}
 Observe that the function $u$ is not a smooth solution as it is not even differentiable. 
 When $p_i$ is positive, we have a peakon moving to the right, and when $p_i$ is negative the antipeakon moves to the left.  The interesting case appears when there is at least one peakon and one antipeakon, which is the case of wave breaking. Explicit formulas exist, see, e.g., 
 \cite{CH:93,CHH:94,BealsSattingerSzm:99,BealsSattingerSzm:00,BealsSattingerSzm:01,parker:08,wahlen}. All these examples are in the so-called conservative case, where the energy is preserved at the wave breaking. However, wave breaking allows for a dichotomy between conservative solutions and dissipative solutions where part of the energy is removed. The analysis of the solution near wave breaking requires a careful change of variables that allows for a smooth transition across wave breaking. For multipeakons this is discussed in \cite{HolRay:06b,HolRay:07B}.  Recently, a new class of solutions was introduced, namely so-called $\alpha$-dissipative solutions that offer a continuous interpolation between conservative ($\alpha=0$) and dissipative ($\alpha=1$) solutions, see \cite{GHR}. In  \cite{GHR}  the symmetric case where the peakon and antipeakon completely annihilate each other at wave breaking, is analyzed in detail.  
 
In this paper we analyze the general case without symmetry, and where the solution does not vanish at wave breaking; in short we extend \cite[Sect.~5]{GHR} to the non-symmetric case. It is somewhat surprising that the non-symmetric case allows an explicit, albeit not simple, solution. The results are presented in this paper.  The crux of the calculation is that one can solve exactly the equation for the characteristics. 

There has also been work on solitary wave solutions of the equation
\begin{equation}\label{eq:CHkappa}
    u_t-u_{txx}+\kappa u_x+3uu_x-2u_xu_{xx}-uu_{xxx}=0,
 \end{equation}
which of course reduces to \eqref{eq:CH} when $\kappa=0$.  The simple transformation $v(t,x)=u(t,x-\kappa t)-\kappa$ takes a solution $u$ of \eqref{eq:CHkappa} into a solution $v$ of \eqref{eq:CH}. If one wants solutions decaying at infinity, this transformation is of little use. However, the decaying solitary wave solutions of \eqref{eq:CHkappa}  do not have the explicit simple form they have for equation \eqref{eq:CH}.  See \cite{cons2001,johnson2003,LiZhang:2004,parker:I,parker:II,parker:III}. A complete description of traveling wave solutions of \eqref{eq:CHkappa} can be found in  \cite{lenells}.

Before we present a summary of the content of this paper, we note that the general Cauchy problem for the Camassa--Holm equation has been extensively studied in both the conservative and dissipative case, see \cite{BreCons:05,BreCons:05a,GHR3,GHR4,GHRb:10,fritz60,GHR5,HolRay:07,HolRay:09}.  In \cite{GHR} the Cauchy problem is studied in the case of a generalized Camassa--Holm system of the form
  \begin{align*}
    u_t-u_{txx}+\kappa u_x+3uu_x-2u_xu_{xx}-uu_{xxx}+\eta\rho\rho_x&=0,\\ 
    \rho_t+(u\rho)_x&=0. 
  \end{align*} 

The present paper is organized as follows. In order to describe the $\alpha$-dissipative peakon-antipeakon solutions, we introduce, in addition to the main unknown $u$, auxiliary variables that measure the concentration of energy. This is done in the form of a Radon measure $\mu$, with the property that its absolutely continuous part satisfies $\muac=u_x^2\,dx$. Whenever wave breaking occurs, $u_x$ tends to $-\infty$ and part of the energy is concentrated on sets of measure zero.  In our case, every solution experiences wave breaking exactly once, and at breaking time, energy is transferred from the absolutely continuous part of $\mu$ to a point mass. By continuing the solution beyond wave breaking without manipulating $\mu$, the peakon and antipeakon are going to pass through each other. However, by removing an $\alpha$-fraction of the energy that is concentrated in a point measure, the solution is either continued by a rescaled peakon-antipeakon for $\alpha\in[0,1)$ or the one-peakon solution for $\
alpha=1$. In addition, we introduce a measure $\nu$ that keeps track of the energy changes. The variables $(u,\mu,\nu)$ are denoted the Eulerian variables. Instead of computing directly the solution in Eulerian coordinates, we are going to introduce Lagrangian variables $(y,U,\bar h,h)$. They are given by the characteristics $y$,  and the Lagrangian velocity defined by 
\begin{equation} 
y_t=u(t,y), \quad U(t,\xi)=u(t,y(t,\xi)),
\end{equation}
as well as two realizations of the energy given by
\begin{equation} 
     \mu=y_\#(\bar h(\xi)\,d\xi),\quad
      \nu=y_\#(h(\xi)\,d\xi).
\end{equation}
There are two main reasons for this change of coordinates. On the one hand  this change of variables associates  the functions $\bar h$ and $h$ to the measures $\mu$ and $\nu$. On the other hand the CH equation rewrites as a system of ordinary differential equations, whose solution remains smooth across wave breaking. At breaking time we modify this system in a continuous manner using the parameter $\alpha$, as introduced in \cite{GHR}. The details of the transformation between the Eulerian and Lagrangian variables can be found in Section \ref{sect:2}, and proofs can be found in \cite{GHR}.  In Section 
\ref{sect:3}  we compute the solution in Eulerian variables, and in Section \ref{sect:4} we provide the detailed calculations in the Lagrangian variables. However, the computation of the full solutions in either set of variables, requires crucial interaction between the two. The fact that the obtained solutions are indeed weak solutions of the CH equation has been established in \cite[Thm.~26]{GHR}, and is not repeated here.

\section{Transformation between Eulerian and Lagrangian variables}\label{sect:2}

The description of $\alpha$-dissipative solutions of the Camassa--Holm equation is based on a generalized method of characteristics. More precisely, one rewrites the CH equation as a system of ordinary differential equations in a suitable Banach space. However, one faces a major problem. Energy can concentrate on sets of measure zero, even in the case of smooth initial data, and hence the corresponding variable is not a function, but rather a measure. To take care of this issue, a suitable change of variables needs to be introduced mapping measures to functions. Thus the aim of this section is  to present the interplay between Eulerian and Lagrangian coordinates.  Since these results are well-established we refer the interested reader to \cite{GHR} for the details and only state the results here. 

We start by introducing the set of Eulerian coordinates $\D$.
\begin{definition}[Eulerian coordinates]
  The set $\D$ is
  composed of all $(u,\mu,\nu)$ such that 
  \begin{enumerate}
  \item $u\in H^1(\Real)$,
  \item $\mu$ is a positive finite
    Radon measure whose absolutely 
    continuous part,  $\muac$,  satisfies 
    \begin{equation}\label{eq:abspart}
      \muac=u_x^2\,dx,
    \end{equation}
  \item $\nu$ is a positive finite Radon measure such that there exists a measurable function $f$ such that 
  \begin{equation}\label{eq:radon}
  \mu = f\nu\quad \text{ and }\quad 0\leq f\leq 1.
\end{equation}
  \end{enumerate}
\end{definition}

Thus any solution of the CH equation is going to be described by a triplet $(u(t,x),\mu(t,x),\nu(t,x))\in \D$, where the measure $\mu$ describes the concentration of energy at breaking times, while the measure $\nu$ needs to be introduced for technical reason, but does not influence the time evolution. 

Let  $G$ be the subgroup of the group of homeomorphisms from $\Real$ to $\Real$ such that 
  \begin{subequations}
    \label{eq:Gcond}
    \begin{align}
      \label{eq:Gcond1}
      f-\id \text{ and } f^{-1}-\id &\text{ both belong to } W^{1,\infty}(\Real), \\
      \label{eq:Gcond2}
      f_\xi-1 &\text{ belongs to } L^2(\Real),
    \end{align}
  \end{subequations}
  where $\id$ denotes the identity function.
  
  Then we can introduce the set of Lagrangian coordinates as follows.
  \begin{definition} \label{def:G} The set $\F$ consists of all $\Theta=(y,U, y_\xi,U_\xi,\bar h, h)$ such that
  \begin{subequations}\label{eq:lagcoord}
    \begin{align}
      \label{eq:lagcoord1}
      &X=(\zeta, U,\zeta_\xi, U_\xi, h)\in L^\infty(\Real)\times [L^2(\Real)\cap L^\infty(\Real)]^4,\\
      \label{eq:lagcoord8}
      & h\in L^1(\Real),\\
      \label{eq:lagcoord3}
      &y_\xi\geq 0, \quad h\geq 0, \quad \bar h\geq 0 \text{  almost everywhere}, \\
      \label{eq:lagcoord4}
      &\lim_{\xi\to -\infty} \zeta(\xi)=0,\\
      \label{eq:lagcoord5}
      &\frac{1}{y_\xi+h}\in L^\infty(\Real),\\ 
      \label{eq:lagcoord6}
      &y_\xi \bar h=U_\xi^2 \text{ almost everywhere},\\ 
      \label{eq:lagcoord7}
      & h\geq \bar h \text{ almost everywhere},\\
      \label{eq:lagcoord9}
      &y+H\in \Gr
    \end{align}
  \end{subequations}
  where we denote $y(\xi)=\zeta(\xi)+\xi$ and $H(t,\xi)=\int_{-\infty}^\xi h(t,\tilde \xi)d\tilde\xi$.
\end{definition}

The condition $y+H\in\Gr$ is crucial since it in general enables to identify equivalence classes and hence enables to identify each element in $\D$ with one equivalence class in $\F$. However, since this will not play a major role for our explicit computations, we will not go into detail here.

The change of variables between Eulerian and Lagrangian coordinates is then given by the following definition. 
\begin{definition} 
  \label{th:Ldef}
  For any $(u,\mu,\nu)$ in $\D$, let
  \begin{subequations}
    \label{eq:Ldef}
    \begin{align}
      \label{eq:Ldef1}
      y(\xi)&=\sup\left\{y\mid \nu((-\infty,y))+y<\xi\right\},\\
      \label{eq:Ldef2}
      h(\xi)&=1-y_\xi(\xi),\\
      \label{eq:Ldef3}
      U(\xi)&=u\circ{y(\xi)},\\
      \label{eq:Ldef5}
      \bar h(\xi)&=f\circ{y(\xi)} h(\xi),
    \end{align}
  \end{subequations}
where $f$ is given through \eqref{eq:radon}.
  Then $\Theta=(y,U,y_\xi, U_\xi, \bar h,h)\in\F$. We denote by $L\colon \D\rightarrow \F$ the mapping which to any
  element $(u,\mu,\nu)\in\D$ associates $\Theta=(y,U,y_\xi, U_\xi,\bar h,h)\in \F$ given by \eqref{eq:Ldef}.
\end{definition}
 For the transformation back to Eulerian variables, we apply the following definition.
  \begin{definition}
  \label{th:umudef} 
  Given any element $\Theta=(y,U,y_\xi, U_\xi,\bar h,h)\in\F$. Then we define $(u,\mu,\nu)$ as follows
  \begin{subequations}
    \label{eq:umudef}
    \begin{align}
      \label{eq:umudef1}
      u(x)&=U(\xi)\text{ for any }\xi\text{ such that  }  x=y(\xi),\\
      \label{eq:umudef2}
      \mu&=y_\#(\bar h(\xi)\,d\xi),\\
      \label{eq:umudef5}
      \nu&=y_\#(h(\xi)\,d\xi).
    \end{align}
  \end{subequations}
  We have that $(u,\mu,\nu)$ belongs to $\D$. We
  denote by $M\colon \F\rightarrow\D$ the mapping
  which to any $\Theta$ in $\F$ associates the element $(u,\mu,\nu)\in \D$ as given
  by \eqref{eq:umudef}.
  \end{definition}
Note that the representation in Lagrangian coordinates depends on which measure $\nu$ we choose. This means, given 
$(u,\mu,\nu)$ and $(u,\mu,\mu)$, then $L(u,\mu,\nu)=(y,U,y_\xi, U_\xi, h,h)\not = L(u,\mu,\mu)=(\tilde y, \tilde U, \tilde y_\xi, \tilde U_\xi , \tilde {\bar {h}}, \tilde h)$. However, if $\supp(\mu_s)=\supp(\nu_s)$,  the support of the singular measures, one can show 
 that there exists a relabeling function $g(\xi)\in G$ such that $\tilde y(\xi)=y(g(\xi))$, $\tilde U(\xi)=U(g(\xi))$ and $\tilde{\bar{h}}(\xi)=h(g(\xi))$, while $\tilde h(\xi)\not =h(g(\xi))$. In addition, this means that the value of $(u,\mu)$ only depends on $(y,U,\bar h)$, but is independent of $h$ in this case.
 
\section{Eulerian coordinates}\label{sect:3}

Consider the following initial data:
\begin{equation}
u(0,x)=p_1(0)e^{-\vert x-q_1(0)\vert}+p_2(0)e^{-\vert x-q_2(0)\vert}, \quad x\in\Real,
\end{equation}
where $p_j(0)$ and $q_j(0)$ are the initial values of the functions
\begin{subequations}
\begin{align}
p_1(t)&=\frac{c_1-c_2e^{L(t-t_0)}}{1-e^{L(t-t_0)}}, \quad p_2(t)=\frac{c_2-c_1e^{L(t-t_0)}}{1-e^{L(t-t_0)}}, \\
q_1(t)&=\ln(L)+c_1(t-t_0)-\ln(c_1-c_2e^{L(t-t_0)}),\\ 
q_2(t)&=-\ln(L)+c_2(t-t_0)+\ln(c_1e^{L(t-t_0)}-c_2).
\end{align}
\end{subequations}
Here $t_0>0$ denotes the future time of wave breaking, which will take place at the origin.  Furthermore, $c_1>0>c_2$, and $L=c_1-c_2$. The fully symmetric case which yields complete annihilation at wave breaking corresponds to $c_1  =-c_2$.\footnote{Most formulas simplify considerably in the fully symmetric case, often after a limiting procedure. It is substantially easier to study the case $c_1 =-c_2$ separately.} Note that $q_1(t)<q_2(t)$ and $q_1(t_0)=q_2(t_0)$.

Then the solution of the CH equation  before wave breaking, which occurs at time $t_0$, is given by  
\begin{equation}\label{sol:Eulerbefore}
u(t,x)=p_1(t)e^{-\vert x-q_1(t)\vert }+p_2(t)e^{-\vert x-q_2(t)\vert}, \quad x\in\Real, \quad t<t_0.
\end{equation}
Define the two Radon measures by 
\begin{equation}
\begin{aligned}
\mu(t)=\nu(t)&=u_x^2(t,x)dx\\
&=\begin{cases} (p_1(t)e^{-q_1(t)}+p_2(t)e^{-q_2(t)})^2e^{2x}dx, &\text{ for }x<q_1(t),\\
(p_2(t)e^{x-q_2(t)}-p_1(t)e^{q_1(t)-x})^2dx, & \text{ for }q_1(t)<x<q_2(t),\\
(p_1(t)e^{q_1(t)}+p_2(t)e^{q_2(t)})^2e^{-2x}dx, &\text{ for }q_2(t)<x.
\end{cases}
\end{aligned}
\end{equation}
Then the energy for $t<t_0$ equals
\begin{equation}\label{eq:2}
\begin{aligned}
\int_\Real \big(u^2(t,x)+u_x^2(t,x)\big)dx & = \int_{-\infty}^{q_1(t)} 2 (p_1(t)e^{-q_1(t)}+p_2(t)e^{-q_2(t)})^2e^{2x}dx\\
&\quad +\int_{q_1(t)}^{q_2(t)} 2(p_1(t)^2e^{2(q_1(t)-x)}+p_2(t)^2e^{2(x-q_2(t))})dx\\
&\quad +\int_{q_2(t)}^{\infty} 2(p_1(t)e^{q_1(t)}+p_2(t)e^{q_2(t)})^2e^{-2x}dx\\
& =2c_1^2+2c_2^2=E^2.
\end{aligned}
\end{equation}
Next we compute the delicate behavior at breaking time when  $t=t_0$.  The solution 
$u$ looks like a one peakon solution with height $c_1+c_2$, that is, 
\begin{equation}
u(t_0,x)=(c_1+c_2)e^{-\vert x\vert }
\end{equation}
as $q_1(t_0)=q_2(t_0)=0$ and $p_j(t)\to\pm\infty$ when $t\to\infty$.  The special case when $c_1+c_2=0$ yields of course the trivial solution. 

Furthermore, we find that for any $M\subset\R$ measurable that
\begin{equation}
\begin{aligned}
\mu(t)(M)=\nu(t)(M)& =\int_M u_x^2(t,x)dx \to (c_1+c_2)^2 \int_M e^{-2\vert x\vert} dx  -4c_1c_2\delta_0(t_0)(M)
\end{aligned}
\end{equation}
where $\delta_0(t_0)(M)=1$ if $0\in M$ and zero otherwise.  The part of the energy which has not concentrated at the origin is given by 
\begin{equation}\label{eq:1}
\begin{aligned}
\int_\Real \big(u^2(t_0,x)+u_x^2(t_0,x)\big)dx&=2\int_\Real (c_1+c_2)^2e^{-2\vert x\vert }dx\\
& = 2(c_1+c_2)^2=E^2+4c_1c_2.
\end{aligned}
\end{equation}
Note that $4c_1c_2$ is negative and hence, comparing \eqref{eq:2} and \eqref{eq:1}, yields that the amount of energy concentrated at the origin at time $t_0$ equals $-4c_1c_2$. 

Thus at $t=t_0$ we find
\begin{equation}
\begin{aligned}
u(t_0,x)&=(c_1+c_2)e^{-\vert x\vert}, \\
\mu(t_0-)=\lim_{t\uparrow t_0}\mu(t_0)&=(c_1+c_2)^2e^{-2\vert x\vert}dx- 4c_1c_2\delta_0(t_0).
\end{aligned}
\end{equation}

Next we introduce the parameter $\alpha\in[0,1]$ that describes $\alpha$-dissipative solutions. In the case of $\alpha=0$, we have the conservative solution where all energy is preserved, while $\alpha=1$ corresponds to the fully dissipative case where all energy concentrated at the origin is removed. 
Then an $\alpha$-fraction of the energy concentrated at the origin at time $t_0$ is given by $-4\alpha c_1c_2$. We modify the energy concentrated at the origin in the measure $\mu$ while keeping it unchanged in $\nu$. More precisely,
\begin{align}
\nu(t_0)&=u_x^2(t_0,x)dx-4c_1c_2\delta_0(t_0), \\
\mu(t_0)&=u_x^2(t_0,x)dx-4(1-\alpha)c_1c_2\delta_0(t_0).
\end{align}

To continue the solution for $t>t_0$ we use as initial data 
\begin{equation}
u(t_0,x)=(c_1+c_2)e^{-\vert x\vert}, \quad
\mu(t_0)=(c_1+c_2)^2e^{-2\vert x\vert}dx- 4c_1c_2\delta_0(t_0).
\end{equation}
We discuss the fully dissipative case separately, and commence with the general case.
\subsection*{The case $\alpha\in[0,1)$:}
The solution will again be a two-peakon solution of the form
\begin{equation}\label{breaksol}
u(t,x)=\tilde p_1(t)e^{-\vert x-\tilde q_1(t)\vert}+\tilde p_2(t)e^{-\vert x-\tilde q_2(t)\vert}, \quad x\in\Real,\quad t>t_0,
\end{equation}
where 
\begin{subequations}
\begin{align}
\tilde p_1(t)&=\frac{d_2-d_1e^{-\tilde L(t-t_0)}}{1-e^{-\tilde L(t-t_0)}},\quad \tilde p_2(t)=\frac{d_1-d_2e^{-\tilde L(t-t_0)}}{1-e^{-\tilde L(t-t_0)}} ,\\
\tilde q_1(t)&=\ln(\tilde L)+d_2(t-t_0)-\ln(d_1e^{-\tilde L(t-t_0)}-d_2),\\
\tilde q_2(t)&=-\ln(\tilde L)+d_1(t-t_0)+\ln(d_1-d_2e^{-\tilde L(t-t_0)}).
\end{align}
\end{subequations}
It remains to compute the values of $d_i$ in terms of $c_i$.  Here  $\tilde L=d_1-d_2$, and it will turn out that $d_1>0>d_2$. In particular, we have $\tilde q_1(t)\leq \tilde q_2(t)$ for all $t\geq t_0$. Furthermore, the energy for $t>t_0$ equals
\begin{equation}
\int_\Real \big(u^2(t,x)+u_x^2(t,x)\big)dx =2d_1^2+2d_2^2.
\end{equation}

We are now ready to establish the connection between the pairs $(c_1,c_2)$ and $(d_1,d_2)$.
By construction we must have 
\begin{equation}
(c_1+c_2)e^{-\vert x\vert}=\lim_{t\to t_0-}u(t,x)=\lim_{t\to t_0+} u(t,x)=(d_1+d_2)e^{-\vert x \vert},
\end{equation}
which implies that 
\begin{equation}
d_1+d_2=c_1+c_2.
\end{equation}
Moreover, we have for the energy, since we take out an $\alpha$-fraction of the energy concentrated at the origin at time $t_0$, that 
\begin{equation}
\begin{aligned}
2d_1^2+2d_2^2& =\lim_{t\to t_0+}\int_\Real (u^2(t,x)+u_x^2(t,x))dx\\& =\lim_{t\to t_0-}\int_\Real (u^2(t,x)+u_x^2(t,x))dx+4\alpha c_1c_2=2c_1^2+2c_2^2+4\alpha c_1c_2.
\end{aligned}
\end{equation}
Thus, $d_1$ satisfies the following quadratic equation,
\begin{equation}
d_1^2-d_1(c_1+c_2)+(1-\alpha)c_1c_2=0.
\end{equation}
Since we assume that $d_1>0$, we get 
\begin{equation}\label{eq:d1d2}
\begin{aligned}
d_1&=\frac{1}{2}(c_1+c_2)+\sqrt{\frac{1}{4}(c_1+c_2)^2-(1-\alpha) c_1c_2} \, , \\
d_2&=\frac{1}{2}(c_1+c_2)-\sqrt{\frac{1}{4}(c_1+c_2)^2-(1-\alpha) c_1c_2}\, .
\end{aligned}
\end{equation}
Note that $d_2<0$, since $-(1-\alpha)c_1c_2>0$.  Furthermore, observe that in the conservative case with $\alpha=0$, we have $d_j=c_j$. The new energy is given by
\begin{equation}
\tilde E^2=2(d_1^2+d_2^2)=E^2+4\alpha c_1c_2.
\end{equation}

Up to this point we are only able to write down the Radon measure $\mu(t,x)$ for $t>t_0$, 
\begin{equation*}
\begin{aligned}
\mu(t,x)=u_x^2(t,x)= \begin{cases}
(\tilde p_1(t)e^{-\tilde q_1(t)}+\tilde p_2(t)e^{-\tilde q_2(t)})^2e^{2x}dx, &\text{ for }x<\tilde q_1(t),\\
(\tilde p_2(t)e^{x-\tilde q_2(t)}-\tilde p_1(t)e^{\tilde q_1(t)-x})^2dx, & \text{ for }\tilde q_1(t)<x<\tilde q_2(t),\\
(\tilde p_1(t)e^{\tilde q_1(t)}+\tilde p_2(t)e^{\tilde q_2(t)})^2e^{-2x}dx, &\text{ for }\tilde q_2(t)<x.
\end{cases}
\end{aligned}
\end{equation*}
The other measure $\nu(t,x)$ is best computed from the solution in Lagrangian coordinates. However, for the sake of completeness, we state the result 
already here (cf.~\eqref{eq:nuBIG}), 
\begin{equation}
\begin{aligned}
\nu(t,x)=\begin{cases} 
\mu(t,x), & \text{ for } x<\tilde q_1(t),\\
\mu(t,x)+\nu_m(t,x), & \text{ for } \tilde q_1(t)<x<\tilde q_2(t),\\
\mu(t,x), & \text{ for } \tilde q_2(t)< x, 
\end{cases}
\end{aligned}
\end{equation}
where 
\begin{equation*}
\begin{aligned}
\nu_m(t,x)&= 4\alpha(1-\alpha)c_1^2c_2^2(1-e^{-\tilde L(t-t_0)})^2e^x\\
&\quad\times\Big(e^{x-d_1(t-t_0)}(\tilde L+d_2e^{d_1(t-t_0)}-d_1e^{d_2(t-t_0)})\\
&\qquad\qquad\qquad\qquad-(d_1-d_2e^{-\tilde L(t-t_0)}-\tilde L e^{d_2(t-t_0)})\Big)^{-2}.
\end{aligned}
\end{equation*}
It should be noted that $\nu_m(t,[\tilde q_1(t), \tilde q_2(t)])=-4\alpha c_1c_2$, and thus equals the amount of energy taken out at time $t=t_0$. Hence, the $\alpha$ part of the energy concentrated at the origin at time $t_0$, is no longer concentrated in one singular point for $t>t_0$, but is, in some sense, spread out over the interval $[\tilde q_1(t), \tilde q_2(t)]$.

\subsection*{The case $\alpha=1$:} In the fully dissipative case $\alpha=1$, the solution reads 
\begin{equation}\label{onepeakon}
u(t,x)=(c_1+c_2)e^{-\vert x-(c_1+c_2)(t-t_0)\vert}, \quad x\in\Real,\quad t>t_0,
\end{equation}
which is the one peakon solution with height and speed equal to $c_1+c_2$. 
Moreover, the energy for $t>t_0$ is given by
\begin{equation}
\int_\Real \big(u^2(t,x)+u_x^2(t,x)\big)dx= 2c_1^2+2c_2^2+4c_1c_2=2(c_1+c_2)^2.
\end{equation}
Here we can write down the measure $\mu(t,x)$, which equals
\begin{equation}
\mu(t,x)=u_x^2(t,x)=(c_1+c_2)^2e^{-2\vert x-(c_1+c_2)(t-t_0)\vert}, \quad x\in\Real, \quad t>t_0.
\end{equation}
As far as the second measure $\nu(t,x)$ is concerned, we again have to postpone the necessary computations until the next section (see \eqref{eq:alpha0stor3}), but state the result here,
\begin{equation}
\nu(t,x)=\mu(t,x)-4c_1c_2\delta_{(c_1+c_2)(t-t_0)}(t), \quad x\in\Real, \quad t>t_0.
\end{equation}

\section{Lagrangian coordinates}\label{sect:4}

Next we turn to the Lagrangian variables, which for $t<t_0$ are solutions of the following system of ordinary differential equations,
\begin{subequations}\label{eq:syslagcoord}
\begin{align}
y_t&=U,\\
U_t&= -Q,\\
y_{t,\xi}&=U_\xi,\\
U_{t,\xi}& =\frac12 h+(U^2-P)y_\xi,\\
h_t&= 2(U^2-P)U_\xi,\\
\bar h_t&=h_t,
\end{align}
\end{subequations}
where $P$ and $Q$ are given by 
\begin{equation}\label{eq:P}
P(t,\xi)=\frac 14 \int_\Real e^{-\vert y(t,\xi)-y(t,\eta)\vert }(2U^2y_\xi+h)(t,\eta)d\eta,
\end{equation}
and 
\begin{equation}\label{eq:Q}
Q(t,\xi)=-\frac14 \int_\Real \sign(\xi-\eta)e^{-\vert y(t,\xi)-y(t,\eta)\vert}(2U^2y_\xi+h)(t,\eta)d\eta.
\end{equation}
The function $\bar h$ is an auxiliary variable whose meaning will only become clear after wave breaking.
Note that before wave breaking, $\bar h$ and $h$ coincide, that is, 
\begin{equation}
\bar h(t,\xi)=h(t,\xi), \quad \text{ for all }\xi\in\Real \text{ and }t<t_0.
\end{equation}
However, the above system is difficult to solve directly, even in the case of the symmetric peakon-antipeakon solution. Instead of solving \eqref{eq:syslagcoord} directly, we will determine the solution by using the connection between Eulerian and Lagrangian variables directly. The key relations are 
\begin{equation}
y_t=u\circ y, \quad U=u\circ y, \quad h=u_x^2\circ y y_\xi.
\end{equation}
We have to determine the initial characteristic (here denoted by $\bar y_0$), given by 
\begin{equation}
\bar y_0(\xi)=\sup\{ y\mid \nu_0((-\infty,y))+y<\xi\}.
\end{equation}
In our example, the measure $\nu_0$ is absolutely continuous and hence the initial characteristic is implicitly given by 
\begin{equation}\label{eq:charac}
\int_{-\infty}^{\bar y_0(\xi)} u_x^2(0,x)dx +\bar y_0(\xi)=\xi.
\end{equation}
Unfortunately, equation \eqref{eq:charac} is hard to solve for $\bar y_0(\xi)$.  However, its derivative is straightforward
\begin{equation}
\begin{aligned}
y_0'(\xi)& =\frac{1}{1+u_{0,x}^2(\bar y_0(\xi))}\\
& =\begin{cases}
(1+(p_1(0)e^{-q_1(0)}+p_2(0)e^{-q_2(0)})^2e^{2\bar y_0(\xi)})^{-1}, &\text{ for }\xi<\xi_1,\\
(1+(p_2(0)e^{\bar y_0(\xi)-q_2(0)}-p_1(0)e^{q_1(0)-\bar y_0(\xi)})^2)^{-1}, & \text{ for }\xi_1<\xi<\xi_2,\\
(1+(p_1(0)e^{q_1(0)}+p_2(0)e^{q_2(0)})^2e^{-2\bar y_0(\xi)})^{-1}, &\text{ for }\xi_2<\xi
\end{cases}\\
& =\begin{cases}
(1+\frac{1}{L^2}(c_1^2e^{c_1t_0}-c_2^2e^{c_2t_0})^2e^{2\bar y_0(\xi)})^{-1}, &\text{ for }\xi<\xi_1,\\
(1-\frac{L^2}{(1-e^{-Lt_0})^2}(e^{c_2t_0+\bar y_0(\xi)}+e^{-(c_1t_0+\bar y_0(\xi))})^2)^{-1},& \text{ for }\xi_1<\xi<\xi_2,\\
 (1+\frac{1}{L^2}(c_1^2e^{-c_1t_0}-c_2^2 e^{-c_2t_0})^2e^{-2\bar y_0(\xi)})^{-1}, & \text{ for }\xi_2<\xi,
\end{cases}
\end{aligned}
\end{equation}
where we introduced $\xi_1$ and $\xi_2$ as  the solutions of $\bar y_0(\xi_1)=q_1(0)$ and $\bar y_0(\xi_2)=q_2(0)$, respectively. For reasons that will become clear later, we will benefit from having characteristics that satisfy $y_0(q_1(0))=q_1(0)$ and $y_0(q_2(0))=q_2(0)$, which is not automatically satisfied by \eqref{eq:charac}. We use the freedom given to us by relabeling to modify $\bar y_0(\xi)$. To that end define 
\begin{equation}
f(z)=\int_{-\infty}^z u_x^2(0,x)dx +z.
\end{equation}
Then $f(z)$ is a relabeling function. Observe that with this definition $\xi_1=f(q_1(0))$, $\xi_2=f(q_2(0))$ and $f'(z)=u_x^2(0,z)+1$. Introduce
\begin{equation}
y_0(z)=\bar y_0(f(z)),
\end{equation}
which implies 
\begin{equation}
y_0(q_1(0))=\bar y_0(f(q_1(0))=\bar y_0(\xi_1)=q_1(0),
\end{equation}
and 
\begin{equation}
y_0(q_2(0))=\bar y_0(f(q_2(0))=\bar y_0(\xi_2)=q_2(0).
\end{equation}
Hence,
\begin{equation}
y_0'(\xi)=\bar y_0'\circ f(\xi)f'(\xi)=1.
\end{equation}
Thus the relabeled initial characteristic is simply $y_0(\xi)=\xi$. Clearly, we could have chosen this function immediately, and the above argument shows that one can always use the identity as the initial characteristic when the initial data contains no singular part. However, the above argument illustrates the actual use of relabeling. 

The key observation here is that the equation for characteristics, $y_t=u\circ y$, can be solved explicitly. Thus,  the Lagrangian variables are then given, using \eqref{eq:syslagcoord} for $t<t_0$, by 
\begin{align*}
y(t,\xi)&= \begin{cases}
y_l(t,\xi), & \text{ for }\xi<q_1(0),\\
y_m(t,\xi), & \text{ for }q_1(0)<\xi<q_2(0),\\
y_r(t,\xi), &\text{ for }q_2(0)<\xi, 
\end{cases}\\
y_{\xi}(t,\xi)& = \begin{cases}
y_{\xi,l}(t,\xi), & \text{ for }\xi<q_1(0),\\
y_{\xi,m}(t,\xi), & \text{ for }q_1(0)<\xi<q_2(0),\\
y_{\xi,r}(t,\xi), &\text{ for }q_2(0)<\xi, 
\end{cases}\\
U(t,\xi)& =\begin{cases}
U_l(t,\xi), & \text{ for } \xi<q_1(0), \\
U_m(t,\xi),  & \text{ for }q_1(0)<\xi<q_2(0)\\
U_r(t,\xi), &\text{ for }q_2(0)<\xi,
\end{cases}\\
h(t,\xi)& =\begin{cases}
h_l(t,\xi), & \text{ for } \xi<q_1(0),\\
h_m(t,\xi), & \text{ for } q_1(0)<\xi<q_2(0),\\
h_r(t,\xi), & \text{ for } q_2(0)<\xi,
\end{cases}
\end{align*}
where 
\begin{align*}
y_l(t,\xi)&= 
\xi+\ln(L)\\
&\quad -\ln\left(L+(c_1e^{-c_1(t-t_0)}-c_1e^{c_1t_0}-c_2e^{-c_2(t-t_0)}+c_2e^{c_2t_0})e^{\xi}\right),\\
y_m(t,\xi)&=\ln\left(\frac{e^{c_2(t-t_0)}}{L}\frac{(c_1e^{L(t-t_0)}-c_2)D(\xi)+L^2e^{c_1(t-t_0)}C(\xi)}{D(\xi)+(c_1e^{c_2(t-t_0)}-c_2e^{c_1(t-t_0)})C(\xi)}\right),\\
y_r(t,\xi)&=\xi-\ln(L)\\
&\quad+\ln\left(L+(c_1e^{c_1(t-t_0)}-c_1e^{-c_1t_0}-c_2e^{c_2(t-t_0)}+c_2e^{-c_2t_0})e^{-\xi}\right), \\
\intertext{and}
y_{\xi,l}(t,\xi)& =
\Big(L+(c_1e^{-c_1(t-t_0)}-c_1e^{c_1t_0}-c_2e^{-c_2(t-t_0)}+c_2e^{c_2t_0})e^\xi\Big)^{-1}L, \\
y_{\xi,m}(t,\xi)& =c_1^2c_2^2 Le^{c_2(t-t_0)}(1-e^{-Lt_0})^2(1-e^{L(t-t_0)})^2 e^\xi\\
&\quad\times\big(D(\xi)+(c_1e^{c_2(t-t_0)}-c_2e^{c_1(t-t_0)})C(\xi)\big)^{-1}\\
&\quad\times\big((c_1e^{L(t-t_0)}-c_2)D(\xi)+L^2e^{c_1(t-t_0)}C(\xi)\big)^{-1}, \\
y_{\xi,r}(t,\xi)& =\Big(L+(c_1e^{c_1(t-t_0)}-c_1e^{-c_1t_0}-c_2e^{c_2(t-t_0)}+c_2e^{-c_2t_0})e^{-\xi}\Big)^{-1}L, \\
\intertext{and}
U_l(t,\xi)& =
\frac{c_1^2e^{-c_1(t-t_0)}-c_2^2e^{-c_2(t-t_0)}}{L+(c_1e^{-c_1(t-t_0)}-c_1e^{c_1t_0}-c_2e^{-c_2(t-t_0)}+c_2e^{c_2t_0})e^\xi}e^\xi,  \\
U_m(t,\xi)& =\Big(D(\xi)^2(c_1^2e^{L(t-t_0)}-c_2^2)+2C(\xi)D(\xi)L^2e^{c_1(t-t_0)}(c_1+c_2)\\
&\qquad +C(\xi)^2L^2e^{c_1(t-t_0)}(c_1^2e^{c_2(t-t_0)}-c_2^2e^{c_1(t-t_0)})\Big)\\
&\quad\times\big((c_1e^{L(t-t_0)}-c_2)D(\xi)+L^2e^{c_1(t-t_0)}C(\xi)\big)^{-1}\\
&\quad\times\big(D(\xi)+(c_1e^{c_2(t-t_0)}-c_2e^{c_1(t-t_0)})C(\xi)\big)^{-1}, \\
U_r(t,\xi)& =\frac{c_1^2e^{c_1(t-t_0)}-c_2^2e^{c_2(t-t_0)}}{L+(c_1e^{c_1(t-t_0)}-c_1e^{-c_1t_0}-c_2e^{c_2(t-t_0)}+c_2e^{-c_2t_0})e^{-\xi}}e^{-\xi}, \\
\intertext{and}
h_l(t,\xi)& =
\frac{(c_1^2e^{-c_1(t-t_0)}-c_2^2e^{-c_2(t-t_0)})^2Le^{2\xi}}{(L+(c_1e^{-c_1(t-t_0)}-c_1e^{c_1t_0}-c_2e^{-c_2(t-t_0)}+c_2e^{c_2t_0})e^\xi)^3}, \\
h_m(t,\xi)& =U^2(t,\xi)y_\xi(t,\xi)-4p_1(t)p_2(t)e^{q_1(t)-q_2(t)}y_\xi(t,\xi), \\
h_m(t,\xi)& =\frac{(c_1^2e^{c_1(t-t_0)}-c_2^2e^{c_2(t-t_0)})^2Le^{-2\xi}}{(L+(c_1e^{c_1(t-t_0)}-c_1e^{-c_1t_0}-c_2e^{c_2(t-t_0)}+c_2e^{-c_2t_0})e^{-\xi})^3}.
\end{align*}
Here
\begin{equation}
\begin{aligned}
C(\xi) &= c_2-c_1e^{-Lt_0}+Le^{\xi+c_2t_0},\\
 D(\xi)&=L^2e^{-c_1t_0}-Lc_1e^\xi+Lc_2e^{\xi-Lt_0}.
\end{aligned}
\end{equation}

The above formulas are obtained as follows. One starts by computing the characteristics using \eqref{sol:Eulerbefore}. The cases $\xi<q_1(0)$  and $q_2(0)<\xi$ are more or less straightforward. However, the challenge is to solve $y_t(t,\xi)=u(t,y(t,\xi))$ for $q_1(0)<\xi<q_2(0)$, since one has to rewrite the resulting ordinary differential equation in a suitable way by applying several changes of variables. The equation reads
\begin{equation*}
y_t(t,\xi)=\frac{L}{1-e^{L(t-t_0)}}(e^{c_1(t-t_0)-y(t,\xi)}-e^{-c_2(t-t_0)+y(t,\xi)}).
\end{equation*}
Introducing $x(t,\xi)=e^{y(t,\xi)}$, the equation becomes
\begin{equation*}
x_t(t,\xi)=\frac{L}{1-e^{L(t-t_0)}}(e^{c_1(t-t_0)}-e^{-c_2(t-t_0)}x^2(t,\xi)).
\end{equation*}
By further introducing $v(t,\xi)=-L(1-e^{L(t-t_0)})^{-1}e^{-c_2(t-t_0)}x(t,\xi)$ we find
 \begin{equation*}
v_t(t,\xi)=-R'(t)+R(t)v(t,\xi)+v^2(t,\xi),
\end{equation*}
where
\begin{equation*}
R(t)=\frac{c_1e^{L(t-t_0)}-c_2}{1-e^{L(t-t_0)}}.
\end{equation*}
Defining the logarithmic derivative $v(t,\xi)=-s_t(t,\xi)/s(t,\xi)$, the equation reads
\begin{equation*}
s_{tt}(t,\xi)=(R(t)s(t,\xi))_t, \quad \text{or} \quad s_t(t,\xi)=R(t)s(t,\xi)+I(\xi).
\end{equation*}
Yet another substitution $w(t,\xi)=(1-e^{L(t-t_0)})s(t,\xi)$ turns the equation into
\begin{equation*}
w_t(t,\xi)=-c_2w(t,\xi)+I(\xi)(1-e^{L(t-t_0)})
\end{equation*}
where $I(\xi)$ is a constant of integration. Finally, by defining $z(t,\xi)=e^{c_2(t-t_0)}w(t,\xi)$ we obtain
\begin{equation*}
z_t(t,\xi)=I(\xi)(e^{c_2(t-t_0)}-e^{c_1(t-t_0)}),
\end{equation*}
which can easily be integrated. By returning to the original variables, we find what is denoted $y_m(t,\xi)$. It is not clear at first sight  that the characteristics for $q_1(0)<\xi<q_2(0)$ are well-defined, and one has to check that the argument in the logarithm is positive and bounded. In particular, one can show that $D(\xi)+(c_1e^{c_2(t-t_0)}-c_2e^{c_1(t-t_0)})C(\xi)<0$ and $(c_1e^{L(t-t_0)}-c_2)D(\xi)+L^2e^{c_1(t-t_0)}C(\xi)<0$ for all $\xi\in[q_1(0),q_2(0)]$. 

Once the characteristics $y(t,\xi)$ are known for $t<t_0$, the Lagrangian velocity $U(t,\xi)=u(t,y(t,\xi))$ is straightforward. 

As far as the energy variable $h(t,\xi)$ is concerned, the necessary computations simplify considerably by rewriting the equations. For $\xi<q_1(0)$, one has that $u_x(t, y(t,\xi))=-u(t,y(t,\xi))$, and hence $h(t,\xi)=U^2(t,\xi)y_\xi(t,\xi)$. Similarly, for $q_2(0)<\xi$, one has that $u_x(t,y(t,\xi))=u(t,y(t,\xi))$, and hence $h(t,\xi)=U^2(t,\xi)y_\xi(t,\xi)$. Again the challenging case is $q_1(0)<\xi<q_2(0)$. 
First we calculate the limit as $t\to t_0-$. Thus it suffices for our purposes to observe $u_x^2(t,y(t,\xi))=u^2(t,y(t,\xi))-4p_1(t)p_2(t)e^{q_1(t)-q_2(t)}$ and hence 
$h(t,\xi)=U^2(t,\xi)y_\xi(t,\xi)-4p_1(t)p_2(t)e^{q_1(t)-q_2(t)}y_\xi(t,\xi)$.

The representation we have chosen so far makes it quite easy to compute the limits as $t\to t_0$:
\begin{subequations}
\begin{align}
y(t_0,\xi)& =\begin{cases} 
\begin{array}{l}
\xi+\ln(L)\\-\ln(L+(c_1-c_1e^{c_1t_0}-c_2+c_2e^{c_2t_0})e^\xi)
\end{array}, & \text{ for } \xi<q_1(0),\\
0, & \text{ for }q_1(0)<\xi<q_2(0),\\
\begin{array}{l}
\xi-\ln(L)\\
+\ln(L+(c_1-c_1e^{-c_1t_0}-c_2+c_2e^{-c_2t_0})e^{-\xi})\end{array}, & \text{ for }q_2(0)<\xi,
\end{cases}\\
U(t_0,\xi)& =\begin{cases}
(c_1+c_2)e^{y(t_0,\xi)}, &\text{ for } \xi<q_1(0),\\
c_1+c_2, & \text{ for } q_1(0)<\xi<q_2(0),\\
(c_1+c_2)e^{-y(t_0,\xi)}, & \text{ for }q_2(0)<\xi,
\end{cases}\\
h(t_0,\xi)& =\begin{cases}
\frac{(c_1+c_2)^2L^3e^{2\xi}}{(L+(c_1-c_1e^{c_1t_0}-c_2+c_2e^{c_2t_0})e^\xi)^3}, &\text{ for }\xi<q_1(0),\\[2mm]
\frac{4c_1^2c_2^2(1-e^{-Lt_0})^2e^\xi}{(Le^{-c_1t_0}+c_2-c_1e^{-Lt_0}+(-c_1+c_2e^{-Lt_0}+Le^{c_2t_0})e^\xi)^2}, & \text{ for } q_1(0)<\xi<q_2(0),\\[2mm]
\frac{(c_1+c_2)^2L^3e^{-2\xi}}{(L+(c_1-c_1e^{-c_1t_0}-c_2+c_2e^{-c_2t_0})e^{-\xi})^3}, & \text{ for }q_2(0)<\xi.
\end{cases}
\end{align}
\end{subequations}
 All these limits are obtained by direct computations. As far as $h(t_0,\xi)$ for $\xi\in[q_1(0),q_2(0)]$ is concerned, observe that $y_\xi(t,\xi)\to 0$ for all $\xi\in[q_1(0),q_2(0)]$ as $t\to t_0-$. Thus, 
$\lim_{t\to t_0} h(t,\xi)=-4 \lim_{t\to t_0}  p_1(t)p_2(t)e^{q_1(t)-q_2(t)}y_\xi(t,\xi)$, since $U(t,\xi)\in H^1(\Real)$ for all times, which simplifies the calculations considerably.
Moreover, note that $h(t_0,\xi)=(c_1+c_2)^2e^{-2\vert y(t_0,\xi)\vert}y_\xi(t_0,\xi)$ for $\xi\not\in[q_1(0),q_2(0)]$.

Next we are going to show that also in Lagrangian coordinates we have that the amount of energy concentrated at the origin at time $t=t_0$ equals $-4c_1c_2$. Note that $\{\xi\in\Real\mid y(t_0,\xi)=0\}=[q_1(0),q_2(0)]$, and thus the amount of energy concentrated at the origin is given by the following integral
\begin{align*}
\int_{q_1(0)}^{q_2(0)} h(t_0,\xi)d\xi&=\frac{4c_1^2c_2^2}{-c_1+c_2e^{-Lt_0}+Le^{c_2t_0}}\int_{\frac{c_1c_2(1-e^{-Lt_0})^2}{c_1-c_2e^{-Lt_0}}}^{\frac{c_1c_2e^{-c_2t_0}(1-e^{-Lt_0})^2}{L}} \frac{1}{\eta^2}d\eta\\
& =-4c_1c_2.
\end{align*}

\bigskip
Let us now compute the solution for $t>t_0$: 

\subsection*{The case $\alpha=1$:} The case $\alpha=1$ yields the dissipative solution, thus we will introduce the function $\bar h(t,\xi)$ for $t>t_0$ as follows
\begin{equation}
\bar h(t,\xi)=\begin{cases}
0, & \text{ for }q_1(0)<\xi<q_2(0),\\
h(t,\xi), & \text{ otherwise.}
\end{cases}
\end{equation} 
In addition, we have to redefine our system \eqref{eq:syslagcoord} of ordinary differential equations. To be more explicit, we have to replace $h(t,\xi)$ by $\bar h(t,\xi)$ 
\textit{everywhere} on the right-hand side for all $t\geq t_0$. Note that this also means that $h(t,\xi)$ has to be replaced by $\bar h(t,\xi)$ in \eqref{eq:P} and \eqref{eq:Q}. Observe that even if $h(t,\xi)$ differs from $\bar h(t,\xi)$ only on the interval $[q_1(0),q_2(0)]$, the solution will be influenced for all $\xi\in\Real$, due to $P(t,\xi)$ and $Q(t,\xi)$ on the right-hand side of \eqref{eq:syslagcoord}. 
We now want to conclude that the solution for $t>t_0$ is given by a one-peakon traveling wave  with height $c_1+c_2$. As far as our system of ordinary differential equations is concerned, it reduces, for $\xi\in(q_1(0),q_2(0))$ and $t\geq t_0$, to 
\begin{subequations}\label{sys:middle}
\begin{align}
y_t& =U,\\
U_t&=-Q,\\
y_{t,\xi}& =0,\\
U_{t,\xi}&=0,\\
h_t&=0,\\
\bar h_t&=0.
\end{align}
\end{subequations}
Thus $y_\xi(t,\xi)=0$, $U_\xi(t,\xi)=0$, and $\bar h(t,\xi)=0$ for all $\xi\in [q_1(0), q_2(0)]$ and $t\geq t_0$. In particular, $y(t,\xi)=c(t)$, $U(t,\xi)=c'(t)$, and $Q(t,\xi)=c''(t)$ for all $\xi\in[q_1(0),q_2(0)]$ and $t\geq t_0$, where $c(t)$ denotes some suitable function only depending on $t$. Since both $y(t,\dott)$ and $U(t,\dott)$ are continuous for any $t\geq t_0$, we have 
\begin{equation}\label{onesided}
y(t, q_1(0)-)=y(t,q_2(0)+)\quad \text{and}\quad U(t, q_1(0)-)=U(t,q_2(0)+),  \quad t\geq t_0.
\end{equation}

As far as the solution for $\xi\not\in[q_1(0),q_2(0)]$ is concerned, we are using a squeezing and relabeling argument based on the considerations for $\xi\in[q_1(0), q_2(0)]$ so far. Namely, let for $t\geq t_0$,
\begin{subequations}\label{helphelpfunctions}
\begin{align}
\tilde y(t,\xi)&=\begin{cases} y(t,\xi), & \text{ for }\xi<q_1(0),\\
y(t, \xi+q_2(0)-q_1(0)), & \text{ for }q_1(0)<\xi,
\end{cases}\\
\tilde U(t,\xi)&=\begin{cases} U(t,\xi), & \text{ for }\xi<q_1(0),\\
U(t,\xi+q_2(0)-q_1(0)), & \text{ for } q_1(0)<\xi,
\end{cases}\\
\tilde{\bar{h}}(t,\xi)& =\begin{cases} \bar h(t,\xi), & \text{ for }\xi<q_1(0),\\
\bar h(t,\xi+q_2(0)-q_1(0)),& \text{ for } q_1(0)<\xi,
\end{cases}\\
\tilde h(t,\xi)&=\tilde{\bar {h}}(t,\xi),\\
\tilde P(t,\xi)& =\begin{cases} P(t,\xi), & \text{ for }\xi<q_1(0),\\
P(t,\xi+q_2(0)-q_1(0)), & \text{ for }q_1(0)<\xi,
\end{cases}\\
\tilde Q(t,\xi)& =\begin{cases} Q(t,\xi) & \text{ for }\xi<q_1(0), \\
Q(t,\xi+q_2(0)-q_1(0)), & \text{ for }q_1(0)<\xi,
\end{cases}
\end{align}
\end{subequations}
which means that we have taken out the part of the function where the energy is concentrated. However, due to \eqref{sys:middle} and \eqref{onesided}, both $\tilde y(t,\xi)$, $\tilde U(t,\xi)$, $\tilde P(t,\xi)$, and $\tilde Q(t,\xi)$ are continuous. 
In particular, the triple $(\tilde y(t,\xi), \tilde U(t,\xi), \tilde{\bar{h}}(t,\xi))$ satisfies the following system of ordinary differential equations for $t\geq t_0$, 
\begin{subequations}\label{eq:syslagcoord2}
\begin{align}
\tilde y_t&=\tilde U,\\
\tilde U_t&= -\tilde Q,\\
\tilde y_{t,\xi}&=\tilde U_\xi,\\
\tilde U_{t,\xi}& =\frac12 \tilde{\bar{h}}+(\tilde U^2-\tilde P)\tilde y_\xi,\\
\tilde h_t&= 2(\tilde U^2-\tilde P)\tilde U_\xi,\\
\tilde {\bar{ h}}_t&=\tilde h_t.
\end{align}
\end{subequations}
A close look reveals that the above system coincides with the one describing the conservative solutions of the Camassa--Holm equation. To be sure that the solution of \eqref{eq:syslagcoord2} for $t\geq t_0$ coincides with the one-peakon solution with height $c_1+c_2$, there are two more properties we have to check. On the one hand, it is left to show that 
\begin{equation}
f(\xi)=\tilde y(t_0,\xi)+\int_{-\infty}^\xi \tilde{\bar{h}}(t_0,\eta)d\eta
\end{equation}
is a relabeling function, which guarantees that $(\tilde y(t_0,\xi), \tilde U(t_0,\xi), \tilde{\bar{h}}(t_0,\xi))$ belongs to $\F$, the set of Lagrangian coordinates. On the other hand, if we can also check that $\tilde y(t_0,\xi)$ is a relabeling function, we can map $\tilde y(t_0,\xi)$ to the identity, thereby ensuring that $(\tilde y(t_0,\xi), \tilde U(t_0,\xi), \tilde {\bar{h}}(t_0,\xi))$ is a Lagrangian representation of the one-peakon centered at the origin. Both claims follow from applying \cite[Lemma 3.5]{GHR}. 
This means, in particular, according to \eqref{onepeakon}, that 
\begin{equation}
\tilde U(t,\xi)=u(t,\tilde y(t,\xi))=(c_1+c_2)e^{-\vert \tilde y(t,\xi)-(c_1+c_2)(t-t_0)\vert} \quad \xi\in\Real, \text{ } t>t_0,
\end{equation}
and especially 
\begin{align*}
\tilde Q(t,\xi)& =-\tilde U_t(t,\xi)\\
& = \sgn(\tilde y(t,\xi)-(c_1+c_2)(t-t_0))\tilde U(t,\xi) (\tilde U(t,\xi)-c_1-c_2).
\end{align*}

\begin{figure}\centering
      \includegraphics[width=.45\textwidth]{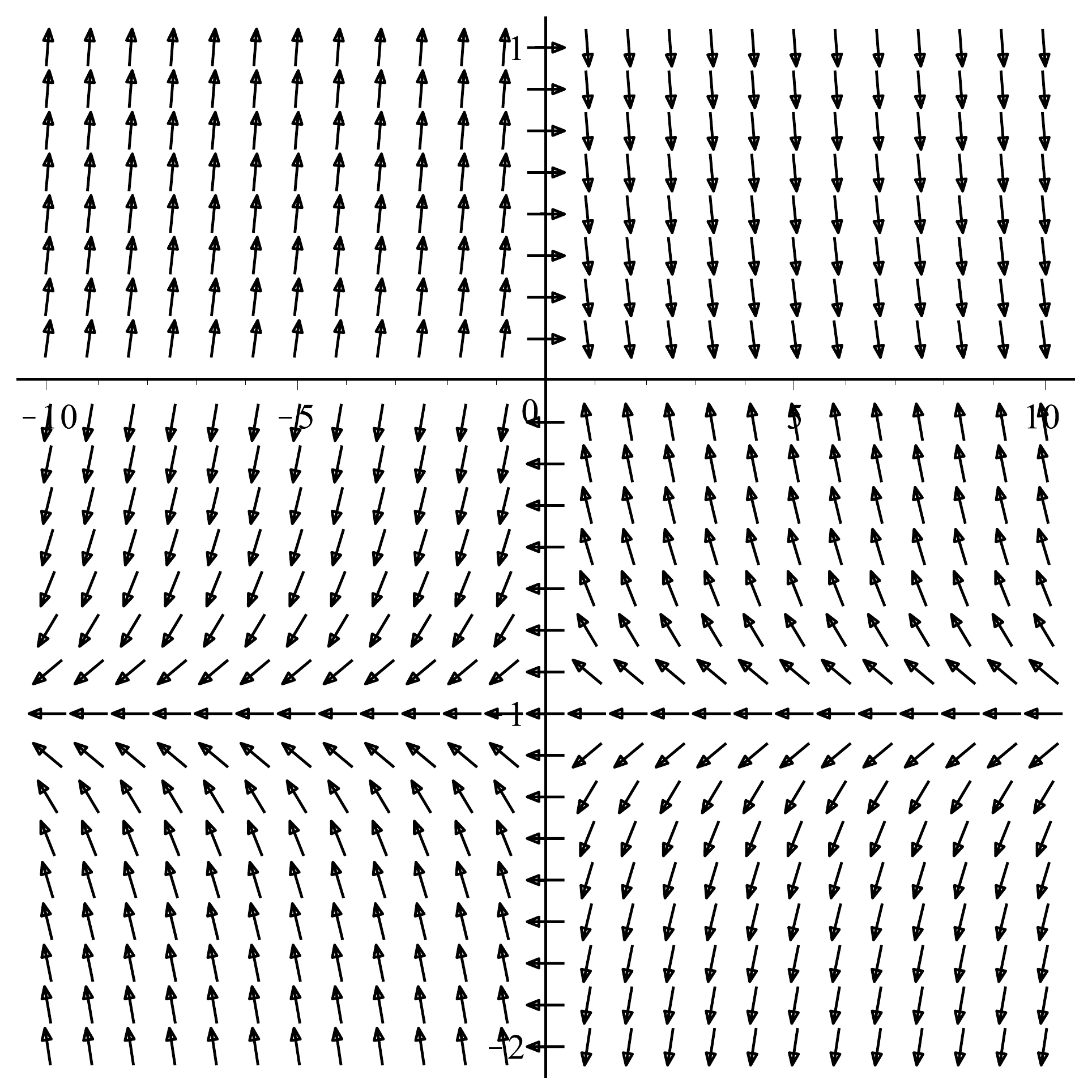} 
      \caption{The vector field for the functions $(z(t),V(t))$ in \eqref{eq:dynsys} for $c_1+c_2=1$.} \label{fig:dynsys}
    \end{figure}

Let $z(t,\xi)=\tilde y(t,\xi)-(c_1+c_2)(t-t_0)$ and $V(t,\xi)=\tilde U(t,\xi)-(c_1+c_2)$.  Then the first two equations in \eqref{eq:syslagcoord2} rewrite as (see Figure \ref{fig:dynsys}) 
\begin{subequations} \label{eq:dynsys}
\begin{align}
z_t& =V,\\
V_t& = -\sgn(z)V(V+(c_1+c_2)),
\end{align}
\end{subequations}
and the above system has a unique solution in $W^{1,\infty}(\Real)\times W^{1,\infty}(\Real)$. This means in particular that we can solve the above system for every $\xi\in\Real$. 

Given $\xi\in\Real$ such that $(z(t_0,\xi), V(t_0,\xi))\not=(0,0)$, then $(z(t,\xi),V(t,\xi))\not =(0,0)$ for any finite time $t\geq t_0$. In particular, if $(z(t_0,\xi), V(t_0,\xi))=(0,0)$ for some $\xi\in\Real$, then $(z(t,\xi),V(t,\xi))=(0,0)$ for all $t\geq t_0$. This means, in particular, that the peakon is traveling along a characteristic, i.e., 
\begin{align}
\tilde y(t,q_1(0))& = (c_1+c_2)(t-t_0),\\
\tilde U(t,q_1(0))& =(c_1+c_2).
\end{align}
Furthermore, we have for all $\xi\in[q_1(0), q_2(0)]$, due to \eqref{onesided} and \eqref{helphelpfunctions}, that
\begin{align}
y(t,\xi)&=y(t, q_1(0)-)=\tilde y(t,q_1(0))=(c_1+c_2)(t-t_0),\\
U(t,\xi)&=U(t, q_1(0)-)=\tilde U(t, q_1(0))=(c_1+c_2),\\
Q(t,\xi)&=Q(t, q_1(0)-)=\tilde Q(t, q_1(0))=0.
\end{align}

As a byproduct of our analysis we  obtained a lot more information about our solution in Eulerian coordinates. Indeed, we have for $t>t_0$,
\begin{subequations} \label{eq:alpha0stor}
\begin{align}
u(t,x)&=(c_1+c_2)e^{-\vert x-(c_1+c_2)(t-t_0)\vert} , \label{eq:alpha0stor1} \\
\mu(t,x)& =u_x^2(t,x)dx,  \label{eq:alpha0stor2}\\
\nu(t,x)& =u_x^2(t,x)dx-4c_1c_2\delta_{(c_1+c_2)(t-t_0)}(t), \label{eq:alpha0stor3}
\end{align}
\end{subequations}
since the peak is traveling at speed $c_1+c_2$. 

The solution in Lagrangian coordinates for $t>t_0$ reads, 
\begin{align}
y(t,\xi)&=\begin{cases}
\begin{array}{l}\xi+\ln(L)\\
-\ln(L+(Le^{-(c_1+c_2)(t-t_0)}-c_1e^{c_1t_0}+c_2e^{c_2t_0})e^\xi),
\end{array} & \text{for $\xi<q_1(0)$},\\
(c_1+c_2)(t-t_0),&\text{for $q_1(0)<\xi<q_2(0)$}, \\
\begin{array}{l}
\xi-\ln(L)\\+\ln(L+(Le^{(c_1+c_2)(t-t_0)}-c_1e^{-c_1t_0}+c_2e^{-c_2t_0})e^-\xi),
\end{array} &\text{for $q_2(0)<\xi$},
\end{cases}\\
U(t,\xi)&=\begin{cases}
\frac{(c_1+c_2)Le^{-(c_1+c_2)(t-t_0)} e^\xi}{L+(Le^{-(c_1+c_2)(t-t_0)}-c_1e^{c_1t_0}+c_2e^{c_2t_0})e^\xi}, &\text{ for }\xi<q_1(0),\\[2mm]
(c_1+c_2), & \text { for }q_1(0)<\xi<q_2(0),\\[2mm]
\frac{(c_1+c_2)Le^{(c_1+c_2)(t-t_0)} e^{-\xi}}{L+(Le^{(c_1+c_2)(t-t_0)}-c_1e^{-c_1t_0}+c_2e^{-c_2t_0})e^{-\xi}}, & \text{ for }q_2(0)<\xi,
\end{cases}\\
h(t,\xi)& =\begin{cases}
\frac{(c_1+c_2)^2L^3e^{-2(c_1+c_2)(t-t_0)}e^{2\xi}}{(L+(Le^{-(c_1+c_2)(t-t_0)}-c_1e^{c_1t_0}+c_2e^{c_2t_0})e^\xi)^3}, &\text{ for } \xi<q_1(0),\\[2mm]
\frac{4c_1^2c_2^2(1-e^{-Lt_0})^2e^\xi}{(Le^{-c_1t_0}+c_2-c_1e^{-Lt_0}+(-c_1+c_2e^{-Lt_0}+Le^{c_2t_0})e^\xi)^2}, & \text{ for } q_1(0)<\xi<q_2(0),\\[2mm]
\frac{(c_1+c_2)^2L^3e^{2(c_1+c_2)(t-t_0)}e^{-2\xi}}{(L+(Le^{(c_1+c_2)(t-t_0)}-c_1e^{-c_1t_0}+c_2e^{-c_2t_0})e^{-\xi})^3}, & \text{ for }q_2(0)<\xi,
\end{cases}\\
\bar h(t,\xi)& =\begin{cases}
\frac{(c_1+c_2)^2L^3e^{-2(c_1+c_2)(t-t_0)}e^{2\xi}}{(L+(Le^{-(c_1+c_2)(t-t_0)}-c_1e^{c_1t_0}+c_2e^{c_2t_0})e^\xi)^3}, &\text{ for } \xi<q_1(0),\\[2mm]
0, & \text{ for } q_1(0)<\xi<q_2(0),\\[2mm]
\frac{(c_1+c_2)^2L^3e^{2(c_1+c_2)(t-t_0)}e^{-2\xi}}{(L+(Le^{(c_1+c_2)(t-t_0)}-c_1e^{-c_1t_0}+c_2e^{-c_2t_0})e^{-\xi})^3}, & \text{ for }q_2(0)<\xi.
\end{cases}
\end{align}

\bigskip
Next we study the general case.
\subsection*{The case $\alpha\in[0,1)$:} 
Here we have to introduce the function $\bar h(t,\xi)$ for $t>t_0$, as follows
\begin{equation}
\bar h(t,\xi)=\begin{cases} h(t,\xi)-\alpha h(t_0,\xi), & \text{ for } q_1(0)<\xi<q_2(0),\\
 h(t,\xi), & \text{ otherwise.}
 \end{cases}
 \end{equation}
 Observe that in the fully conservative case, $\alpha=0$, we have $\bar h(t,\xi)=h(t,\xi)$.
We have to redefine our system  \eqref{eq:syslagcoord} of ordinary differential equations for $t\geq t_0$. To be more explicit, we have to replace $h(t,\xi)$  by $\bar h(t,\xi)$ \textit{everywhere} on the right-hand side for all $t\geq t_0$. Note that this also means that  the function $h(t,\xi)$ has to be replaced by $\bar h(t,\xi)$ in the definitions of $P$ and $Q$, i.e.,   in \eqref{eq:P} and \eqref{eq:Q}, respectively. Note that for $\alpha\in(0,1)$, $h(t,\xi)$ differs from $\bar h(t,\xi)$ only on the interval $[q_1(0),q_2(0)]$.  However, due to the nonlocal nature of $P(t,\xi)$ and $Q(t,\xi)$ on the right-hand side of \eqref{eq:syslagcoord}, the solution will be influenced for all $\xi\in\Real$!

Again it is difficult to solve \eqref{eq:syslagcoord} with $h(t,\xi)$ replaced by $\bar h(t,\xi)$ on the right-hand side for $t\geq t_0$. Thus we proceed as follows.  We know that the solution in Eulerian coordinates at breaking time is given by 
\begin{subequations}\label{initial:alphadiss}
\begin{align}
u(t_0,x)& =(d_1+d_2)e^{-\vert x\vert}=(c_1+c_2)e^{-\vert x\vert},\\
\mu(t_0,x)& =u_x^2(t_0,x) dx-4d_1d_2\delta_0(t_0)=u_x^2(t_0,x)dx-4(1-\alpha)c_1c_2\delta_0(t_0),\\
\nu(t_0,x) & = u_x^2(t_0,x) dx -4c_1c_2\delta_0(t_0).
\end{align}
\end{subequations}

\vspace{0.2cm}
In the special case $\alpha=0$, we obtain the well-studied conservative solution, that is, $d_j=c_j$ for $j=1$, $2$ and the solution for $t>t_0$ equals \eqref{breaksol} with $\tilde p_j=p_j$ and $\tilde q_j=q_j$ for $j=1$, $2$. In other words, the solution \eqref{sol:Eulerbefore} is valid for all $t\geq t_0$. Moreover, in this case $\bar h(t,\xi)=h(t,\xi)$ for all $t\geq t_0$ and hence $\mu(t)=\nu(t)$ for all $t\geq t_0$. Thus the measure $\nu$ is not needed as it does not add any information. 

\vspace{0.2cm}
For $\alpha\in (0,1)$ the situation is a bit more involved. Let $(y(t,\xi), U(t,\xi), \bar h(t,\xi), h(t,\xi))$ be the $\alpha$-dissipative solution in Lagrangian coordinates. Since we replaced $h(t,\xi)$ by $\bar h(t,\xi)$ everywhere on the right-hand side of \eqref{eq:syslagcoord} for $t\geq t_0$, we have that the time evolution of $(y(t,\xi), U(t,\xi), \bar h(t,\xi), h(t,\xi))$ is independent of $h(t,\xi)$. Thus if 
\begin{equation}
y(t_0,\xi)+\int_{-\infty}^\xi \bar h(t_0,\eta) d\eta= y(t_0,\xi)+\bar H(t_0,\xi) \in G,
\end{equation}
then $(y(t,\xi), U(t,\xi), \bar h(t,\xi), \bar h(t,\xi))$ is the solution of \eqref{eq:syslagcoord} with (valid!) initial data $(y(t_0,\xi), U(t_0,\xi), \bar h(t_0,\xi), \bar h(t_0,\xi))$ for $t\geq t_0$. However, $(y(t,\xi), U(t,\xi), \bar h(t,\xi), \bar h(t,\xi))$ in Lagrangian coordinates corresponds to the conservative solution in Eulerian coordinates with initial data $(u(t_0), \mu(t_0), \mu(t_0))$, which is given according to the case $\alpha=0$, by \eqref{sol:Eulerbefore} with $c_j$ replaced by $d_j$ for $j=1$, $2$. 

 Thus it is left to show that 
$y(t_0,\xi)+\bar H(t_0,\xi)$ is a relabeling function. We apply the fundamental lemma \cite[Lemma 3.5]{GHR}, which reduces this difficult task to showing that there exists $c>0$ (which may depend on $t_0$) such that 
\begin{equation}\label{eq:test}
c<y_\xi(t_0,\xi)+\bar h(t_0,\xi),\quad \text{ for all } \xi\in\Real.
\end{equation}
To that end we observe that
\begin{equation}
y_\xi(t_0,\xi)\geq \begin{cases}
\frac{1}{L}(c_1e^{c_1t_0}-c_2e^{c_2t_0}), & \text{ for } \xi<q_1(0),\\
0, & \text{ for } q_1(0)<\xi<q_2(0),\\
\frac{1}{L} (c_1e^{-c_1t_0}-c_2e^{-c_2t_0}), & \text{ for } q_2(0)<\xi,
\end{cases}
\end{equation}
and 
\begin{equation}
\bar h(t_0,\xi)\geq \begin{cases}
0, & \text{ for } \xi<q_1(0),\\
4d_1d_2e^{q_1(0)}\min(c_1-c_2e^{-Lt_0}, e^{c_2t_0} L), & \text{ for } q_1(0)<\xi<q_2(0),\\
0, & \text{ for } q_2(0)<\xi.
\end{cases}
\end{equation}
Thus
\begin{equation*}
y_\xi(t_0,\xi)+\bar h(t_0,\xi)\geq \begin{cases} \frac{1}{L} (c_1e^{c_1t_0}-c_2e^{c_2t_0}), & \text{ for } \xi<q_1(0),\\
 4d_1d_2e^{q_1(0)}\min(c_1-c_2e^{-Lt_0}, e^{c_2t_0} L), & \text{ for } q_1(0)<\xi<q_2(0),\\
 \frac{1}{L} (c_1e^{-c_1t_0}-c_2e^{-c_2t_0}), & \text{ for } q_2(0)<\xi,
 \end{cases}
 \end{equation*}
 which proves \eqref{eq:test}, since all terms on the right-hand side are strictly positive due to our assumption that $c_1>0>c_2$ and hence also $d_1>0>d_2$.
 
Finally recall that $h(t,\xi)=\bar h(t,\xi)+\alpha h(t_0,\xi)$, which at first enables us to compute $h(t,\xi)$ and in a further step to derive $\nu(t)$.

We find that the solution in Lagrangian coordinates reads
\begin{align*}
y(t,\xi)&=\begin{cases}
y_l(t,\xi), & \text{ for } \xi<q_1(0),\\
y_m(t,\xi), & \text{ for } q_1(0)<\xi<q_2(0),\\
y_r(t,\xi), & \text{ for } q_2(0)<\xi,
\end{cases}\\
y_\xi(t,\xi) & =\begin{cases}
y_{\xi,l}(t,\xi), & \text { for } \xi<q_1(0),\\
y_{\xi,m}(t,\xi), & \text{ for } q_1(0)<\xi<q_2(0),\\
y_{\xi,r}(t,\xi), &\text{ for } q_2(0)<\xi,
\end{cases}\\
U(t,\xi)&=\begin{cases}
U_l(t,\xi), & \text{ for } \xi<q_1(0),\\
U_m(t,\xi), & \text{ for }q_1(0)<\xi<q_2(0),\\
U_r(t,\xi), & \text{ for } q_2(0)<\xi,
\end{cases}\\
\bar h(t,\xi)&=\begin{cases}
\bar h_l(t,\xi), & \text{ for } \xi<q_1(0),\\
\bar h_m(t,\xi), & \text{ for }q_1(0)<\xi<q_2(0),\\
\bar h_r(t,\xi), & \text{ for } q_2(0)<\xi,
\end{cases}\\
h(t,\xi)&=\begin{cases}
h_l(t,\xi), & \text{ for } \xi<q_1(0),\\
h_m(t,\xi), & \text{ for }q_1(0)<\xi<q_2(0),\\
h_r(t,\xi), & \text{ for } q_2(0)<\xi.
\end{cases}
\end{align*}
Here
\begin{align*}
y_l(t,\xi)&=\xi +\ln(L\tilde L)\\
&\quad-\ln(L\tilde L+(\tilde L(-c_1e^{c_1t_0}+c_2e^{c_2t_0})+L(d_1e^{-d_1(t-t_0)}-d_2e^{-d_2(t-t_0)}))e^\xi),\\
y_m(t,\xi)&=\ln\Big(-e^{d_1(t-t_0)}\frac{\tilde Le^{d_2(t-t_0)}(S(\xi)-1)+(-d_1+d_2e^{-\tilde L(t-t_0)})(S(\xi)+1)}{\tilde L(S(\xi)+1)+(d_2e^{d_1(t-t_0)}-d_1e^{d_2(t-t_0)})(S(\xi)-1)}\Big),\\
y_r(t,\xi)&=\xi-\ln(L\tilde L)\\
&\quad+\ln(L\tilde L+(\tilde L(-c_1e^{-c_1t_0}+c_2e^{-c_2t_0})+L(d_1e^{d_1(t-t_0)}-d_2e^{d_2(t-t_0)}))e^{-\xi}), \\
\intertext{and}
y_{\xi,l}(t,\xi)&=
\frac{L\tilde L}{L\tilde L+(\tilde L(-c_1e^{c_1t_0}+c_2e^{c_2t_0})+L(d_1e^{-d_1(t-t_0)}-d_2e^{-d_2(t-t_0)}))e^\xi)},\\
y_{\xi,m}(t,\xi)&=-2d_1d_2e^{d_1(t-t_0)}(1-e^{-\tilde L(t-t_0)})^2 S'(\xi)\\
&\quad\times\Big(\tilde L(S(\xi)+1)+(d_2e^{d_1(t-t_0)}-d_1e^{d_2(t-t_0)})(S(\xi)-1)\Big)^{-1}\\
&\quad\times\Big((d_1-d_2e^{-\tilde L(t-t_0)})(S(\xi)+1)-\tilde L e^{d_2(t-t_0)}(S(\xi)-1)\Big)^{-1}, \\
y_{\xi,r}(t,\xi)&=\frac{L\tilde L}{L\tilde L+(\tilde L(-c_1e^{-c_1t_0}+c_2e^{-c_2t_0})+L(d_1e^{d_1(t-t_0)}-d_2e^{d_2(t-t_0)}))e^{-\xi}}, \\
\intertext{and}
U_l(t,\xi)&=
\frac{(d_1^2e^{-d_1(t-t_0)}-d_2^2e^{-d_2(t-t_0)})Le^\xi}{L\tilde L+(\tilde L(-c_1e^{c_1t_0}+c_2e^{c_2t_0})+L(d_1e^{-d_1(t-t_0)}-d_2e^{-d_2(t-t_0)}))e^\xi}, \\
U_m(t,\xi)&=\tilde L \Big((d_1^2-d_2^2e^{-\tilde L(t-t_0)})(S(\xi)+1)^2-2(d_1^2-d_2^2)e^{d_2(t-t_0)}(S(\xi)^2-1)\\
&\qquad\qquad\qquad\qquad+e^{d_2(t-t_0)}(d_1^2e^{d_2(t-t_0)}-d_2^2e^{d_1(t-t_0)})(S(\xi)-1)^2\Big)\\
&\quad\times\Big((\tilde L(S(\xi)+1)+(d_2e^{d_1(t-t_0)}-d_1e^{d_2(t-t_0)})(S(\xi)-1)\Big)^{-1}\\
&\quad\times\Big((d_1-d_2e^{-\tilde L(t-t_0)})(S(\xi)+1)-\tilde Le^{d_2(t-t_0)}(S(\xi)-1))\Big)^{-1}, \\
U_r(t,\xi)&=\frac{(d_1^2e^{d_1(t-t_0)}-d_2^2e^{d_2(t-t_0)})Le^{-\xi}}{L\tilde L+(\tilde L (-c_1e^{-c_1t_0}+c_2e^{-c_2t_0})+L(d_1e^{d_1(t-t_0)}-d_2e^{d_2(t-t_0)}))e^{-\xi}}, \\
\intertext{and}
\bar h_l(t,\xi)&=
\frac{(d_1^2e^{-d_1(t-t_0)}-d_2^2e^{-d_2(t-t_0)})^2L^3\tilde L e^{2\xi}}{(L\tilde L+(\tilde L(-c_1e^{c_1t_0}+c_2e^{c_2t_0})+L(d_1e^{-d_1(t-t_0)}-d_2e^{-d_2(t-t_0)}))e^\xi)^3},\\
\bar h_m(t,\xi)&=U^2(t,\xi)y_\xi(t,\xi)-4\tilde p_1(t)\tilde p_2(t)e^{\tilde q_1(t)-\tilde q_2(t)}y_\xi(t,\xi),\\
\bar h_r(t,\xi)&=\frac{(d_1^2e^{d_1(t-t_0)}-d_2^2e^{d_2(t-t_0)})^2L^3\tilde Le^{-2\xi}}{(L\tilde L+(\tilde L (-c_1e^{-c_1t_0}+c_2e^{-c_2t_0})+L(d_1e^{d_1(t-t_0)}-d_2e^{d_2(t-t_0)}))e^{-\xi})^3},\\
\intertext{and}
h_l(t,\xi)&=
\frac{(d_1^2e^{-d_1(t-t_0)}-d_2^2e^{-d_2(t-t_0)})^2L^3\tilde L e^{2\xi}}{(L\tilde L+(\tilde L(-c_1e^{c_1t_0}+c_2e^{c_2t_0})+L(d_1e^{-d_1(t-t_0)}-d_2e^{-d_2(t-t_0)}))e^\xi)^3}, \\
h_m(t,\xi)&=U^2(t,\xi)y_\xi(t,\xi)-4\tilde p_1(t)\tilde p_2(t)e^{\tilde q_1(t)-\tilde q_2(t)}y_\xi(t,\xi)-2\alpha c_1c_2S'(\xi), \\
h_r(t,\xi)&=\frac{(d_1^2e^{d_1(t-t_0)}-d_2^2e^{d_2(t-t_0)})^2L^3\tilde Le^{-2\xi}}{(L\tilde L+(\tilde L (-c_1e^{-c_1t_0}+c_2e^{-c_2t_0})+L(d_1e^{d_1(t-t_0)}-d_2e^{d_2(t-t_0)}))e^{-\xi})^3},
\end{align*}
where 
\begin{align*}
Q(t_0+,\xi)& =(1-\alpha)Q(t_0-,\xi)\\
&= d_1d_2\Big(\frac{2c_1c_2L(1-e^{-Lt_0})^2-(c_1-c_2e^{-Lt_0}+Le^{c_2t_0})(LC(\xi)+D(\xi))}{(-c_1+c_2e^{-Lt_0}+Le^{c_2t_0})(LC(\xi)+D(\xi))}\Big)
\end{align*}
and 
\begin{equation*}
S(\xi)=\frac{Q(t_0+,\xi)}{d_1d_2}=\frac{2c_1c_2L(1-e^{-Lt_0})^2-(c_1-c_2e^{-Lt_0}+Le^{c_2t_0})(LC(\xi)+D(\xi))}{(-c_1+c_2e^{-Lt_0}+Le^{c_2t_0})(LC(\xi)+D(\xi))}.
\end{equation*}
Note that $S(\xi)\in[-1,1]$ for $\xi\in[q_1(0),q_2(0)]$.
Moreover, direct computations yield
\begin{equation}
S'(\xi)=-2c_1c_2\frac{L^2(1-e^{-Lt_0})^2 e^\xi}{(LC(\xi)+D(\xi))^2}.
\end{equation}

\begin{figure}\centering
      \includegraphics[width=.45\textwidth]{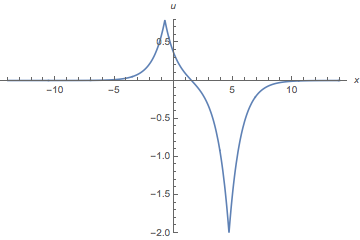} \\
       \includegraphics[width=.45\textwidth]{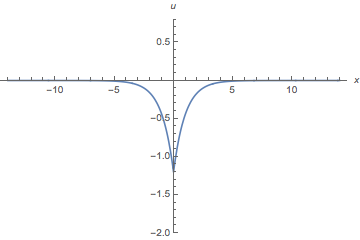} \\
      \includegraphics[width=.45\textwidth]{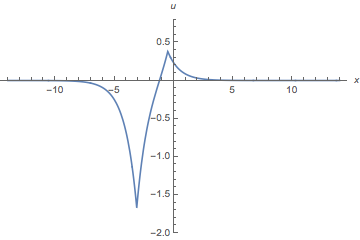}   \\
      \includegraphics[width=.45\textwidth]{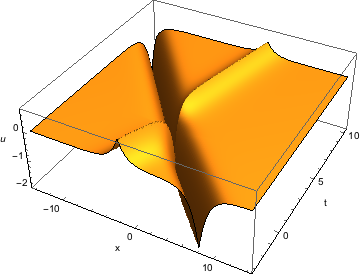}  
      \caption{The function $u$ is plotted for times $t=-1.5$, $t=1.0=t_0$, and $t=3.0$. Parameter values $c_1=0.8$, $c_2=-2.0$, $\alpha=0.5$. } \label{fig:u}
    \end{figure}

\begin{figure}\centering
      \includegraphics[width=.45\textwidth]{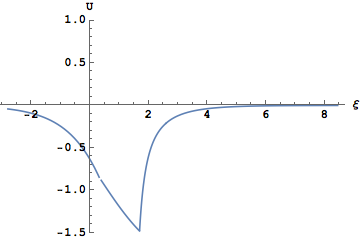} \\
       \includegraphics[width=.45\textwidth]{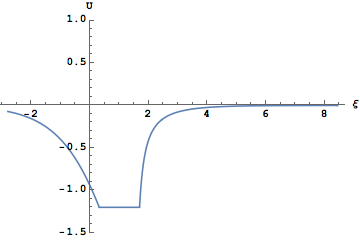} \\
      \includegraphics[width=.45\textwidth]{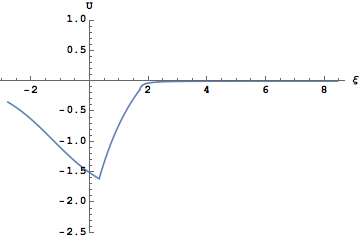}   \\
      \includegraphics[width=.45\textwidth]{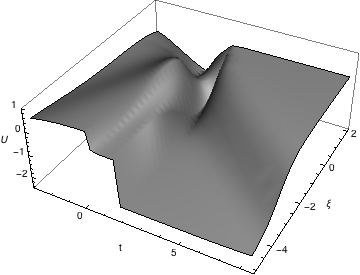}  
      \caption{The function $U$ is plotted for times $t=-0.8$, $t=1.0=t_0$, and $t=2.0$. Parameter values as in Fig. \ref{fig:u}. } \label{fig:U}
    \end{figure}

\begin{figure}\centering
      \includegraphics[width=.85\textwidth]{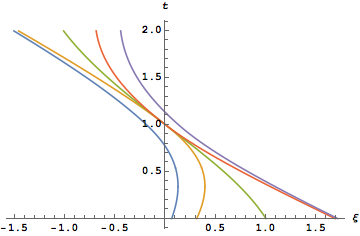} 
      \caption{The characteristics are plotted for five different values of $\xi$. Parameter values as in Fig. \ref{fig:u}. } \label{fig:y}
    \end{figure}

\begin{figure}\centering
      \includegraphics[width=.45\textwidth]{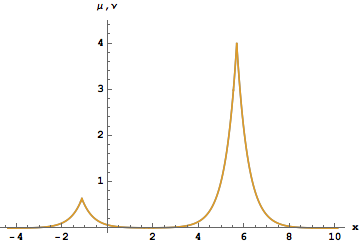} \\
      \includegraphics[width=.45\textwidth]{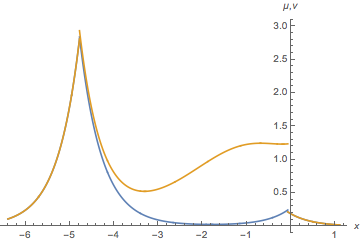}
      \caption{The measures $\mu,\nu$ are plotted for $t=-3.0$ (they coincide) and $t=4.0$. The measure $\nu$ is discontinuous at $x=\tilde q_j(t)$ for $j=1,2$ and $t>t_0$. 
      Parameter values as in Fig. \ref{fig:u}. } \label{fig:mu}
    \end{figure}

In this case the solution in Eulerian coordinates reads for $t>t_0$ (recall \eqref{breaksol})
\begin{align}
u(t,x)&=\tilde p_1(t)e^{-\vert x-\tilde q_1(t)\vert}+\tilde p_2(t) e^{-\vert x-\tilde q_2(t)\vert }, \notag\\
\mu(t,x)&=u_x^2(t,x) dx, \notag\\
\nu(t,x)&=\begin{cases}
u_x^2(t,x) dx , & \text{ for } x<\tilde q_1(t),\\
u_x^2(t,x)dx+\nu_m(t,x) dx, & \text{ for } \tilde q_1(t)<x<\tilde q_2(t),  \label{eq:nuBIG}\\
u_x^2(t,x) dx , & \text{ for } \tilde q_2(t)<x,
\end{cases}
\end{align}
where
\begin{equation*}
\begin{aligned}
\nu_m(t,x)&= 4\alpha(1-\alpha)c_1^2c_2^2(1-e^{-\tilde L(t-t_0)})^2e^x\\
&\quad\times\Big(e^{x-d_1(t-t_0)}(\tilde L+d_2e^{d_1(t-t_0)}-d_1e^{d_2(t-t_0)})\\
&\qquad\qquad\qquad\qquad\qquad-(d_1-d_2e^{-\tilde L(t-t_0)}-\tilde L e^{d_2(t-t_0)})\Big)^{-2}.
\end{aligned}
\end{equation*}
The solution $u$ is displayed in Figure \ref{fig:u}, while $U$ is plotted in Figure \ref{fig:U}. The characteristics  $y$ are visualized in Figure \ref{fig:y}, and the measures can be found in Figure \ref{fig:mu}.

\section*{Acknowledgement}  Thanks are due to Hans Lundmark who convinced us this was worthwhile.
%
%


\begin{thebibliography}{99}


 \bibitem{BealsSattingerSzm:99}
     \newblock R.~Beals, D.~Sattinger, and J.~Szmigielski. 
     \newblock Multi-peakons and a theorem of Stieltjes. 
     \newblock {\em Inverse Problems}  15:1--4, 1999. 

\bibitem{BealsSattingerSzm:00}
     \newblock R.~Beals, D.~Sattinger, and J.~Szmigielski. 
     \newblock Multipeakons and the classical moment problem. 
     \newblock {\em Adv. Math.}  154:229--257, 2000. 



   \bibitem{BealsSattingerSzm:01}
     \newblock R.~Beals, D.~Sattinger, and J.~Szmigielski. 
     \newblock Peakon-antipeakon interaction. 
     \newblock {\em J. Nonlinear Math. Phys.}  8:23--27, 2001. 

  
\bibitem{BreCons:05}
  \newblock A.~Bressan and A.~Constantin. 
  \newblock Global conservative solutions of the Camassa--Holm equation. 
  \newblock {\em Arch. Ration. Mech. Anal.},  183:215--239, 2007. 

\bibitem{BreCons:05a}
  \newblock A.~Bressan and A.~Constantin. 
  \newblock Global dissipative solutions of the Camassa--Holm equation. 
  \newblock {\em Analysis and Applications},  5:1--27, 2007. 
  
\bibitem{CH:93}
  \newblock R.~Camassa and D.~D. Holm. 
  \newblock An integrable shallow water equation with peaked solitons. 
  \newblock {\em Phys. Rev. Lett.}, 71(11):1661--1664, 1993. 

\bibitem{CHH:94}
\newblock R.~Camassa, D.~D. Holm, and J. M. Hyman.
\newblock A new integrable shallow water equation. 
\newblock {\em Adv. Appl. Mech.} 31:1--33, 1994.


   \bibitem{cons2001}
     \newblock A.~Constantin. 
     \newblock On the scattering problem for the Camassa--Holm equation. 
     \newblock {\em Proc. R. Soc. A } 457:953--970, 2001.


   \bibitem{cons_esc1}
     \newblock A.~Constantin and J.~Escher. 
     \newblock Wave breaking for nonlinear nonlocal shallow water equations. 
     \newblock {\em  Acta Math.}, 181:229--243, 1998. 


\bibitem{consstrauss}
   \newblock A.~Constantin and W.~Strauss. 
\newblock Stability of peakons.
\newblock {\em Comm. Pure Appl. Math.}  53:603--610, 2000.
 
   \bibitem{GHR3}
 \newblock     K. Grunert, H. Holden, and X. Raynaud.
     \newblock Global conservative solutions of the
  Camassa--Holm equation for initial data  with  nonvanishing asymptotics.
   \newblock   {\em Discrete Cont. Dyn. Syst., Series A}, 32:4209--4277, 2012. 

\bibitem{GHR4}
 \newblock K. Grunert, H. Holden, and X. Raynaud.
  \newblock Global solutions for the two-component Camassa--Holm system.
  \newblock  {\em Comm. Partial Differential Equations}, 37:2245--2271, 2012. 

\bibitem{GHRb:10} K. Grunert, H. Holden, and  X. Raynaud. 
  \newblock Lipschitz metric for the  {C}amassa--{H}olm equation on the  line. 
  \newblock  {\em Discrete Contin. Dyn. Syst.} 33:2809--2827, 2013.

\bibitem{fritz60}
 \newblock K. Grunert, H. Holden, and X. Raynaud.
\newblock Periodic conservative solutions for the two-component Camassa--Holm system.
\newblock  In {\em  Spectral Analysis, Differential Equations and Mathematical Physics.} 
 {\em A Festschrift for Fritz Gesztesy on the Occasion of his 60th Birthday} (eds. H. Holden, B. Simon, and G. Teschl)
Amer. Math. Soc., pp.~165--182, 2013. 


\bibitem{GHR5}
  \newblock K. Grunert, H. Holden, and X. Raynaud.
  \newblock Global dissipative solutions of the two-component Camassa--Holm system for initial data with nonvanishing asymptotics. 
  \newblock {\em Nonlinear Anal. Real World Appl.} 17:203--244, 2014.

\bibitem{GHR}
 \newblock K. Grunert, H. Holden, and X. Raynaud.
\newblock A continuous interpolation between conservative and dissipative solutions for the two-component Camassa--Holm system.
\newblock {\em Forum of Mathematics, Sigma}  vol. 1, e1, 70 pages doi:10.111, 2014.


\bibitem{HolRey:2006}
   \newblock H.~Holden and X.~Raynaud. 
\newblock A convergent numerical  scheme for the
Camassa--Holm equation based on multipeakon
\newblock {\em Discrete and Continuous Dynamical System} 14:505--523, 2006.   
   
   \bibitem{HolRay:06b}
    \newblock H.~Holden and X.~Raynaud. 
     \newblock Global conservative multipeakon solutions of the Camassa--Holm equation. 
     \newblock {\em J. Hyperbolic Differ. Equ.},  4:39--64, 2007. 

 
\bibitem{HolRay:07}
  \newblock H. Holden and X. Raynaud. 
  \newblock Global conservative solutions of the {C}amassa--{H}olm equation --- a
  Lagrangian point of view. 
  \newblock {\em Comm. Partial Differential Equations},  32:1511--1549, 2007. 

\bibitem{HolRay:07B}
\newblock H.~Holden and X.~Raynaud. 
\newblock Global dissipative multipeakon solutions for the Camassa--Holm equation 
\newblock {\em Commun. in Partial Differential Equations}, 33:2040--2063, 2008.

\bibitem{HolRay:09}
  \newblock H. Holden and X. Raynaud. 
  \newblock Dissipative solutions of the {C}amassa--{H}olm equation. 
  \newblock {\em Discrete Cont. Dyn. Syst.},  24:1047--1112, 2009. 


 \bibitem{johnson2003}
     \newblock R. S. Johnson. 
     \newblock On solutions of the Camassa--Holm equation. 
     \newblock {\em Proc. R. Soc. A } 459:1687--1708, 2003.


\bibitem{lenells}
\newblock J. Lenells.
\newblock Traveling wave solutions of the Camassa--Holm equation.
\newblock {\em J. Differential Eqn.} 217:393--430, 2005.

\bibitem{LiZhang:2004}
\newblock Y. Li and J. E. Zhang.
\newblock The multiple-soliton solution of the Camassa--Holm equation.
  \newblock {\em Proc. R. Soc. A } 460:2617--2627, 2004.
  

\bibitem{parker:I}
\newblock A. Parker.
\newblock On the Camassa--Holm equation and a direct method of solution. I. Bilinear form and solitary waves.
\newblock {\em Proc. R. Soc. A } 460:2929--2957, 2004.

\bibitem{parker:II}
\newblock A. Parker.
\newblock On the Camassa--Holm equation and a direct method of solution. II. Soliton solutions.
\newblock {\em Proc. R. Soc. A } 461:3611--3632, 2005.


\bibitem{parker:III}
\newblock A. Parker.
\newblock On the Camassa--Holm equation and a direct method of solution. III. $N$-soliton solutions.
\newblock {\em Proc. R. Soc. A } 461:3893--3911, 2005.

\bibitem{parker:08}
\newblock A. Parker.
\newblock Wave dynamics for peaked solitons of the Camassa--Holm equation.
\newblock {\em Chaos, Solitons and Fractals} 35:220--237, 2008.

  
 \bibitem{wahlen}
 \newblock E. Wahl\'en.
\newblock The interaction of peakons and antipeakons.
\newblock {\em Dyn.  Contin. Discrete Impuls. Syst. Ser. A.} 13:465-472, 2006.





\end{thebibliography}
\end{document}